\input amstex
\documentstyle{amsppt}
\magnification=\magstep1                        
\hsize6.5truein\vsize8.9truein                  
\NoRunningHeads
\loadeusm

\magnification=\magstep1                        
\hsize6.5truein\vsize8.9truein                  
\NoRunningHeads
\loadeusm

\document
\topmatter

\title  
Coppersmith-Rivlin type inequalities and the order of vanishing of polynomials at $1$
\endtitle

\rightheadtext{the multiplicity of the zero at $1$ of polynomials}

\author Tam\'as Erd\'elyi
\endauthor


\address Department of Mathematics, Texas A\&M University,
College Station, Texas 77843, College Station, Texas 77843 (T. Erd\'elyi) \endaddress


\thanks {{\it 2010 Mathematics Subject Classifications.} 11C08, 41A17, 26C10, 30C15}
\endthanks

\keywords 
\endkeywords

\abstract
For $n \in {\Bbb N}$ , $L > 0$, and $p \geq 1$ let $\kappa_p(n,L)$ be the largest possible value of $k$ for which there is a 
polynomial $P \neq 0$ of the form
$$P(x) = \sum_{j=0}^n{a_jx^j}\,, \qquad |a_0| \geq L \Bigg( \sum_{j=1}^n{|a_j|^p} \Bigg)^{1/p}, \quad a_j \in {\Bbb C}\,,$$
such that $(x-1)^k$ divides $P(x)$. For $n \in  {\Bbb N}$, $L > 0$, and $q \geq 1$ 
$\mu_q(n,L)$ be the smallest value of $k$ for which there is a polynomial $Q$ of degree $k$ with complex
coefficients such that
$$|Q(0)| > \frac 1L \Big( \sum_{j=1}^n{|Q(j)|^q} \Big)^{1/q}\,.$$
We find the size of $\kappa_p(n,L)$ and $\mu_q(n,L)$ for all $n \in {\Bbb N}$, $L > 0$, and $1 \leq p,q \leq \infty$. 
The result about $\mu_\infty(n,L)$ is due to Coppersmith and Rivlin, but our proof is completely different 
and much shorter even in that special case.

\endabstract

\endtopmatter

\document

\head 1. Notation \endhead

In [B-99] and [B-13] we examined a number of problems concerning polynomials with coefficients
restricted in various ways. We are particularly interested in how small such polynomials can be
on the interval $[0,1]$. For example, we proved that there are absolute constants $c_1>0$ and $c_2>0$
such that
$$\exp\left(-c_1\sqrt n \right) \leq \min_{0 \neq Q \in {\Cal F}_n} 
\left\{{\max_{x \in [0,1]}{|Q(x)|}} \right\} \leq \exp\left(-c_2\sqrt n \right)$$
for every $n \geq 2$, where ${\Cal F}_n$ denotes the set of all polynomials of degree at most $n$
with coefficients from $\{-1,0,1\}$.

Littlewood considered minimization problems of this variety on the unit disk.
His most famous, now solved, conjecture was that the $L_1$ norm of an element $f \in {\Cal F}_n$ on
the unit circle grows at least as fast as $c\log N$, where $N$ is the number of non-zero coefficients
in $f$ and $c > 0$ is an absolute constant.

When the coefficients are required to be integers, the questions have a
Diophantine nature and have been studied  from a variety of points of view. See
[A-79, B-98, B-95, F-80, O-93].

One key to the analysis is a study of the related problem of giving an
upper bound for the multiplicity of the zero these restricted polynomials
can have at $1$. In [B-99] and [B-13] we answer this latter question precisely for the class of
polynomials of the form
$$Q(x) = \sum_{j=0}^n {a_jx^j}\,, \qquad |a_j| \leq 1\,, \quad a_j \in {\Bbb C}\,, \quad j=1,2,\ldots, n\,,$$
with fixed $|a_0| \neq 0$.

Variants of these questions have attracted considerable study, though rarely have
precise answers been possible to give. See in particular [A-90, B-32, B-87, E-50, Sch-33, Sz-34].
Indeed, the classical, much studied, and presumably very difficult problem of Prouhet, Tarry,
and Escott rephrases as a question of this variety. (Precisely: what is the maximal vanishing
at $1$ of a polynomial with integer coefficients with $l_1$ norm $2n$? It is conjectured to be
$n$.) See [H-82], [B-94], or [B-02].

For $n \in  {\Bbb N}$, $L > 0$, and $p \geq 1$ we define the following numbers.
Let $\kappa_p(n,L)$ be the largest possible value of $k$ for which there is a polynomial $Q \neq 0$ of the form
$$Q(x) = \sum_{j=0}^n{a_jx^j}\,, \qquad |a_0| \geq L \Bigg( \sum_{j=1}^n{|a_j|^p} \Bigg)^{1/p}, \quad a_j \in {\Bbb C}\,,$$
such that $(x-1)^k$ divides $Q(x)$. For $n \in {\Bbb N}$ and $L > 0$ let $\kappa_\infty(n,L)$ the largest possible value of $k$ 
for which there is a polynomial $Q \neq 0$ of the form
$$Q(x) = \sum_{j=0}^n{a_jx^j}\,, \qquad |a_0| \geq L \max_{1 \leq j \leq n}{|a_j|}\,, \quad a_j \in {\Bbb C}\,,$$
such that $(x-1)^k$ divides $Q(x)$. In [B-13] we proved that there is an absolute constant $c_3 > 0$ such that 
$$\min \Big\{\frac {1}{6} \sqrt{(n(1-\log L)} - 1\,, n \Big\} \leq \kappa_{\infty}(n,L) \leq \min\Big\{c_3\sqrt{n(1-\log L)}\,, n\Big\}$$
for every $n \in {\Bbb N}$ and $L \in (0,1]$. However, we were far from being able to establish the right 
result in the case of $L \geq 1$. In [B-13] we proved the right order of magnitude of   
$\kappa_{\infty}(n,L)$ and $\kappa_2(n,L)$ in the case of $L \geq 1$. 
Our results in [B-99] and [B-13] sharpen and generalize results of Schur [Sch-33],
Amoroso [A-90], Bombieri and Vaaler [B-87], and Hua [H-82] who gave versions of this
result for polynomials with integer coefficients. Our results in [B-99]
have turned out to be related to a number of recent papers from a rather wide range
of research areas.
See [A-02, B-98, B-95, B-96 B-97a, B-97b, B-97, B-00, B-07, B-08a, B-08b, Bu-99, C-02, C-13, C-10,
D-99, D-01, D-03, D-13, E-08a, E-08b, F-00, G-05, K-04, K-09, M-03, M-68, N-94, O-93, P-99, P-12, P-13,
R-04, R-07, S-99, T-07, T-84], for example.
More on the zeros of polynomials with Littlewood-type coefficient constraints
may be found in [E-02b]. Markov and Bernstein type inequalities under Erd\H os type
coefficient constraints are surveyed in [E-02a].

For $n \in  {\Bbb N}$, $L > 0$, and $q \geq 1$ we define the following numbers.
Let $\mu_q(n,L)$ be the smallest value of $k$ for which there is a polynomial of degree $k$ with complex
coefficients such that
$$|Q(0)| > \frac 1L \Big( \sum_{j=1}^n{|Q(j)|^q} \Big)^{1/q}\,.$$
Let $\mu_\infty(n,L)$ be the smallest value of $k$ for which there is a polynomial of degree $k$ with complex
coefficients such that
$$|Q(0)| > \frac 1L \max_{j \in \{1,2,\ldots,n\}}{|Q(j)|}\,.$$
It is a simple consequence of H\"older's inequality (see Lemma 3.6) that
$$\kappa_p(n,L) \leq \mu_q(n,L)\,.$$
whenever $n \in  {\Bbb N}$, $L > 0$, $1 \leq p,q \leq \infty$, and $1/p + 1/q = 1$. 

In this paper we find the the size of $\kappa_p(n,L)$ and $\mu_q(n,L)$ for all $n \in {\Bbb N}$, $L > 0$, and 
$1 \leq p,q \leq \infty$. The result about $\mu_\infty(n.L)$ is due to Coppersmith and Rivlin, [C-92], but our proof 
presented in this paper is completely different and much shorter even in that special case.

\head 2 New Results \endhead

\proclaim{Theorem 2.1}
Let $p \in (1,\infty]$ and $q \in [1,\infty)$ satisfy $1/p + 1/q = 1$. There are absolute constants $c_1 > 0$ and $c_2 > 0$
such that
$$\sqrt{n}(c_1L)^{-q/2} - 1 \leq \kappa_p(n,L) \leq \mu_q(n,L)  \leq \sqrt{n}(c_2L)^{-q/2} + 2$$
for every $n \in {\Bbb N}$ and $L > 1/2$, and
$$c_3 \min \Big\{\sqrt{n(-\log L)}, n \Big\}  \leq \kappa_p(n,L) \leq \mu_q(n,L) \leq c_4 \min \Big\{\sqrt{n(-\log L)},n \Big\} + 4$$
for every $n \in {\Bbb N}$ and $L \in (0,1/2]$. Here $c_1 := 1/53$, $c_2 := 40$, $c_3 := 2/7$, and $c_4 := 13$ are appropriate choices.  
\endproclaim

\proclaim{Theorem 2.2}
There are constants absolute $c_1 > 0$ and $c_2 > 0$ such that
$$ c_1 \sqrt{n(1-L)} - 1 \leq \kappa_1(n,L) \leq \mu_\infty(n,L)  \leq c_2 \sqrt{n(1-L)} + 1$$
for every $n \in {\Bbb N}$ and $L \in (1/2,1]$, and
$$c_3 \min \Big\{ \sqrt{n(-\log L)}, n \Big\} \leq \kappa_1(n,L) \leq \mu_\infty(n,L) \leq c_4 \min \Big\{\sqrt{n(-\log L)},n \Big\} + 4$$
for every $n \in {\Bbb N}$ and $L \in (0,1/2]$.
Note that $\kappa_1(n,L) = \mu_\infty(n,L) = 0$ for every $n \in {\Bbb N}$ and $L > 1$.
Here $c_1 := 1/5$, $c_2 := 1$, $c_3 := 2/7$, and $c_4 := 13$ are appropriate choices.
\endproclaim

\head 3. Lemmas \endhead

In this section we list our lemmas needed in the proofs of Theorems 2.1 and 2.2. These lemmas 
are proved in Section 4. Let ${\Cal P}_n$ be the set of all polynomials of degree at most $n$ with real coefficients. 
Let ${\Cal P}_n^c$ be the set of all polynomials of degree at most $n$ with complex coefficients.

\proclaim{Lemma 3.1} Let $p \in (1,\infty)$. For any $1 \leq M$ there are polynomials $P_n$ of the form
$$P_n(x) = \sum_{j=0}^n{a_{j,n}x^j}\,, \qquad a_{j,n} \in {\Bbb R}\,, \quad a_{0,n} \geq \frac{3M}{\pi^2} + o(M)\,,$$ 
$$\Bigg( \sum_{j=1}^n{|a_{j,n}|^p} \Bigg)^{1/p} \leq 16M^{1/p}\,,$$
such that $P_n$ has at least $\lfloor \sqrt{n/M} \rfloor$ zeros at $1$.
\endproclaim

\proclaim{Lemma 3.2}  Let $p,q \in (1,\infty)$ satisfy $1/p + 1/q = 1$. For any $L \geq 1/48$ there are polynomials $P_n$ of the form
$$P_n(x) = \sum_{j=0}^n{a_{j,n}x^j}\,, \qquad a_{j,n} \in {\Bbb R}\,, \quad a_{0,n} \geq L + o(L) \,, \quad \sum_{j=1}^n{|a_{j,n}|^p} \leq 1\,,$$
such that $P_n$ has at least $\lfloor \sqrt{n}(cL)^{-q/2} \rfloor$ zeros at $1$ with $\displaystyle{c := \frac{3}{16\pi^2}}$.
\endproclaim

\proclaim{Lemma 3.3} Let $p \in [1,\infty)$. For any $L \in (0,1/17)$ there are polynomials $P_n$ of the form
$$P_n(x) = \sum_{j=0}^n{a_{j,n}x^j}\,, \qquad a_{j,n} \in {\Bbb R}\,, \quad a_{0,n} = L\,, \quad \sum_{j=1}^n{|a_{j,n}|^p} \leq 1\,,$$
such that $P_n$ has at least $\displaystyle{\frac 27} \min\{\sqrt{n(1-\log L)},n\}$ zeros at $1$.
\endproclaim

\proclaim{Lemma 3.4} For any $L \in (0,1)$ there are polynomials $P_n \not\equiv 0$ of the form
$$P_n(x) = \sum_{j=0}^n{a_{j,n}x^j}\,, \qquad a_{j,n} \in {\Bbb R}\,, \quad a_{0,n} \geq L \sum_{j=1}^n{|a_{j,n}|}\,,$$
such that $P_n$ has at least $\displaystyle{\frac 15} \sqrt{(n-1)(1-L)}$ zeros at $1$.
\endproclaim

The observation below is well known, easy to prove, and recorded in several papers. See [B-99], for example.

\proclaim{Lemma 3.5}
Let $P \neq  0$ be a polynomial of the form $P(x) = \sum_{j=0}^n{a_jx^j}$.
Then $(x-1)^k$ divides $P$ if and only if $\sum_{j=0}^n{a_jQ(j)} = 0$ for all polynomials
$Q \in {\Cal P}_{k-1}^{c}$.
\endproclaim

Our next lemma is a simple consequence of H\"older's inequality.

\proclaim{Lemma 3.6}
Let $1 \leq p,q \leq \infty$ and $1/p + 1/q = 1$. Then for every $n \in  {\Bbb N}$ and $L > 0$, we have
$$\kappa_p(n,L) \leq \mu_q(n,L)\,.$$
\endproclaim

The next lemma is stated as Lemma 3.4 in [K-03], where a proof of it is also presented.

\proclaim{Lemma 3.7}
For arbitrary real numbers $A,M > 0$, there exists a polynomial $G$ such that $F = G^2 \in {\Cal P}_m$ with 
$$m < \sqrt{\pi} \sqrt{A} \root 4 \of M + 2$$
such that $F(0) = M$ and
$$|F(x)| \leq \min \{ M,x^{-2} \}\,, \qquad x \in (0,A]\,.$$
\endproclaim

We also need Lemma 5.7 from [B-99] which may be stated as follows. 

\proclaim {Lemma 3.8}
Let $n$ and $R$ be positive integers with $1 \leq R \leq \sqrt{n}$.
Then there exists a polynomial $F \in {\Cal P}_m$ with
$$m \leq 4 \sqrt{n} + \textstyle{\frac{9}{7}}R\sqrt{n} + R + 4 \leq 
\textstyle{\frac{44}{7}}R\sqrt{n} + 4$$
such that
$$F(1) = F(2) = \cdots = F(R^2) = 0$$
and
$$|F(0)| > \exp(R^2)\big(|F(R^2+1)| + |F(R^2+2)| + \cdots + |F(n)|\big) \geq \exp(R^2)\Bigg( \sum_{j=1}^n{|F(j)|^2} \Bigg)^{1/2}\,.$$ 
\endproclaim

Lemmas 3.6 and 3.7 imply the following results needed in the proof of Theorems 2.1 and 2.2. 

\proclaim{Lemma 3.9}
Let $q \in [1,\infty)$. For every $n \in {\Bbb N}$, $q \in [1,\infty)$, and $K > 0$, there are polynomials $F \in {\Cal P}_m$ satisfying
$$|F(0)| > K \Bigg( \sum_{j=1}^n{|F(j)|^q} \Bigg)^{1/q} \quad {\text \rm and} \quad 
m \leq \cases \sqrt{n}(40K)^{q/2} + 2\,, \enskip & 0 < K < 2\,, \\
13 \min \Big \{\sqrt{n \log K}, n \Big\} + 4\,, \enskip & K \geq 2\,. 
\endcases$$
\endproclaim

\proclaim{Lemma 3.10}
For every $n \in {\Bbb N}$ and $K > 1$, there are polynomials $F \in {\Cal P}_m$ satisfying
$$|F(0)| > K \max_{j \in \{1,2,\ldots, n\}}{|F(j)|}\,, \quad {\text \rm and} \quad 
m \leq \cases \sqrt{n(K-1)/2} + 1\,, \enskip & 1 < K < 2\,, \\
13 \min \Big \{\sqrt{n \log K}, n \Big\} + 4\,, \enskip & K \geq 2\,. 
\endcases$$

\endproclaim

\head 4. Proofs of the Lemmas \endhead 

\demo{Proof of Lemma 3.1}
Modifying the construction on page 138 of [B-95] we define $H_1(x) := 1$ and 
$$H_m(x) := \frac{(-1)^{m+1}2(m!)^2}{2\pi i} \int_{\Gamma}{\frac{x^t \, dt}{(t-2)\prod_{j=0}^m{(t-j^2)}}}\,, \qquad m=2,3, \ldots \,, \quad x \in (0,\infty)\,,$$
where the simple closed contour $\Gamma$ surrounds the zeros of the denominator of the integrand. Then $H_m$ is 
a polynomial of degree $m^2$ with a zero at $1$ with multiplicity at least $m+1$. (This can be seen easily by repeated 
differentiation and then evaluation of the above contour integral  by expanding the contour to infinity.)  
Also, by the residue theorem,
$$H_m(x) = 1 + d_mx^2 + \sum_{k=1}^m{c_{k,m}x^{k^2}}\,, \qquad m=2,3, \ldots \,, \tag 4.1$$
where
$$c_{k,m} = \frac{(-1)^{m+1}2(m!)^2}{(k^2-2)\prod_{j=0,j \neq k}^m{(k^2-j^2)}} = \frac{4}{k^2-2}\frac{(-1)^{k+1}(m!)^2}{(m-k)!(m+k)!}\,,$$
and
$$d_m = \frac{(-1)^{m+1}2(m!)^2}{\prod_{j=0}^m{(2-j^2)}} \,.$$ 
It follows that each $c_{k,m}$ is real and 
$$|c_{k,m}| \leq \frac{4}{|k^2-2|}\,, \qquad k=1,2,\ldots, m\,, \tag 4.2$$
and a simple calculation shows that 
$$|d_m| \leq 8\,, \qquad  m=2,3, \ldots \,. \tag 4.3$$ 
(No effort has been made to optimize the bound in (4.3).)  
Let $S_M$ be the collection of all odd square free integers in $[1,M]$. Let  $m := \lfloor \sqrt{n/M} \rfloor$. 
If $m = 0$ then there is nothing to prove. So we may assume that $m \geq 1$. 
It is well known that 
$$|S_M| \geq \frac{3M}{\pi^2} + o(M)\,,$$
where $|A|$ denotes the number of elements in a finite set $a$.
This follows from the fact that if $S_M^*$ is the collection of all square free integers in $[1,M]$, then
$$|S_M^*| = \frac{6M}{\pi^2} + o(M)\,,$$ 
see [H-38, pp. 267-268], for example, by observing that the number of odd square free integer in $[1,M]$ is not less than 
the number of even square free integers in $[1,M]$ (if $a$ is an even square free integer then $a/2$ is 
an odd square free integer). We define
$$P_n(x) := \sum_{j \in S_M}{H_m(x^j)}\,.$$ 
Then $P_n$ is of the form 
$$P_n(x) = \sum_{j=0}^n{a_{j,n}x^j}\,, \qquad a_{j,n} \in {\Bbb R}\,, \quad j=0,1,\ldots,n\,.$$
We have
$$a_{0,n} = |S_M^*| \geq \frac{3M}{\pi^2} + o(M)\,.$$
First assume that $m=1$. Then 
$$\sum_{j=1}^n{|a_{j,n}|^p} = 2|S_M| \leq 2M\,,$$  
and as $P_n$ has $1$ zero at $1$, the lemma follows. 
Now assume that $m \geq 2$.
Since $ju \neq lv$ whenever $j,l \in S_M$, $j \neq l$, and $u,v \in \{1^2,2^2,\ldots,m^2\} \cup \{2\}$, we have
$$\split \sum_{j=1}^n{|a_{j,n}|^p} \leq & \, |S_M| \left( 8^p + \sum_{k=1}^m{\Bigg( \frac{4}{|k^2-2|} \Bigg)^p} \right) 
\leq |S_M| \left( 8^p + \sum_{k=1}^m{\frac{4^p}{|k^2-2|}} \right) \cr 
= & \, M(8^p + 8^p) \leq 16^pM\,. \cr \endsplit$$
Observe that each term in $P_n$ has a zero at $1$ with multiplicity at least $m+1 > \lfloor \sqrt{n/M} \rfloor$ zeros at $1$, 
and hence so does $P_n$. 
\qed \enddemo

\demo{Proof of Lemma 3.2}
The statement follows from Lemma 3.1 by choosing $1 \leq M$ so that 
$$L := \frac{3}{16\pi^2}M^{1-1/p} = \frac{3}{16\pi^2}M^{1/q}\,.$$
This can be done when $\displaystyle{\frac{3}{16\pi^2} \leq L}$.
\qed \enddemo

\demo{Proof of Lemma 3.3}
Let $L \in (0,1/17]$. We define
$$k := \min \left\{\left \lfloor \frac{-\log L}{\log 17} \right \rfloor,n \right\} \qquad \text {and} \qquad 
m :=  \lfloor \sqrt{n/k} \rfloor\,.$$  
Observe that $k \geq 1$ and $m \geq 1$ hold.
Let $P_n := LH_m^k \in {\Cal P}_n$, where $H_m \in {\Cal P}_{m^2}$ defined by (4.1). Then
$$P_n(x) = \sum_{j=0}^n{a_{j,n}x^j}\,, \qquad a_{j,n} \in {\Bbb R}\,, \quad j=0,1,\ldots, n\,,$$
has at least 
$$km \geq k \frac 12 \, \sqrt{n/k} = \frac 12 \sqrt{nk} = \frac{1}{2 \sqrt{\log 17}} \min \Big \{\sqrt{n(-\log L)},n \Big\}$$  
zeros at $1$, where $2 \sqrt{\log 17} < 7/2$.
Clearly, $a_{0,n} = P_n(0) = L$, and using the notation in (4.1), we can deduce that
$$\split \sum_{j=1}^n{|a_{j,n}|^p} \leq & L^p \Bigg( \sum_{j=1}^n{|a_{j,n}|} \Bigg)^p \leq L^p \Bigg(1 + |d_m| + \sum_{k=1}^m{|c_{k,m}|} \Bigg)^{kp}  \cr 
\leq & L^p(1 + 8 + 8)^{kp} = L^p \, 17^{kp} \leq L^pL^{-p} = 1\,, \cr \endsplit $$
if $m \geq 2$, and 
$$\sum_{j=1}^n{|a_{j,n}|^p} \leq L^p \Bigg( \sum_{j=1}^n{|a_{j,n}|} \Bigg)^p \leq L^p 2^{kp} \leq L^pL^{-p} = 1\,,$$
if $m=1$.
\qed \enddemo

\demo{Proof of Lemma 3.4}
Let 
$$r := \left \lfloor 12 \, \frac{1+L}{1-L} \right \rfloor + 1\, \qquad \text {and} \qquad m := \left \lfloor \sqrt{\frac{n-1}{r}} \right \rfloor\,.$$
When $m \leq 1$ we have $\lfloor (1/9)\sqrt{n(1-L)} \rfloor = 0$, so there is nothing to prove.
Now assume that $m \geq 2$.
Let $P_n \in {\Cal P}_n$ be defined by $P_n(x) := H_m(x^r)$, where $H_m \in {\Cal P}_{m^2}$ defined by (4.1). 
Let $Q_n \in {\Cal P}_n$ be defined by 
$$Q_n(x) = -\int_0^1{P_n(t) \, dt}  + \int_{0}^x{P_n(t)\,dt}\,.$$
Then, using the notation in (4.1), we have
$$Q_n(x) = -1 - \frac{d_m}{2r+1} - \sum_{k=1}^m{\frac{c_{k,m}}{r k^2+1}} + 
x + \frac{d_mx^{2r +1}}{2r +1} + \sum_{k=1}^m{\frac{c_{k,m}x^{rk^2+1}}{rk^2+1}}\,.$$
Writing 
$$Q_n(x) = \sum_{j=0}^n{a_{j,n}x^j}\,, \qquad a_{j,n} \in {\Bbb R}\,, \quad j=0,1,\ldots, n\,,$$
and recalling (4.2) and (4.3), we have 
$$|a_{0,n}| \geq  1 - \frac{8}{2r +1} - \sum_{k=1}^m{\frac{4}{|k^2-4|(rk^2+1)}} \geq 1- \frac{8}{2r +1} - \frac 8r > 1 - \frac{12}{r}\,,$$
and 
$$\sum_{j=1}^n{|a_{j,n}|} \leq 1 + \frac{8}{2r +1} + \sum_{k=1}^m{\frac{4}{(k^2-2)(rk^2+1)}} < 1 + \frac{12}{r}\,.$$  
Combining the previous two inequalities, we obtain
$$\frac{|a_{0,n}|}{\sum_{j=1}^n{|a_{j,n}|}} > \frac{1-12/r}{1+12/r} \geq \frac{1-(1-L)/(1+L)}{1+(1-L)/(1+L)} = L\,.$$
Also $Q_n$ has at least $m+1 \geq \lfloor \sqrt{(n-1)/r} \rfloor + 1 \geq  \displaystyle{\frac{1}{5}}\sqrt{(n-1)(1-L)}$ zeros at $1$.
\qed \enddemo

\demo{Proof of Lemma 3.6}  We assume that $p,q \in (1,\infty)$, the result in the cases  $p=1, q=\infty$ and $p = \infty, q=1$ can 
be proved similarly with straightforward modification of the proof.
Let $m := \mu_q(n,L)$.  Let $Q$ be a polynomial of degree $m$ with complex coefficients such that 
$$|Q(0)| > \frac 1L \Bigg( \sum_{j=1}^n{|Q(j)|^q} \Bigg)^{1/q}\,.$$
Now let $P$ be a polynomial of the form
$$P(x) = \sum_{j=0}^n{a_jx^j}\,, \qquad |a_0| \geq L \Bigg( \sum_{j=1}^n{|a_j|^p} \Bigg)^{1/p}, \quad a_j \in {\Bbb C}\,.$$
It follows from H\"older's inequality that
$$\Bigg| \sum_{j=1}^n{a_jQ(j)} \Bigg| \leq 
\Bigg( \sum_{j=1}^n{|a_j|^p} \Bigg)^{1/p}\Bigg( \sum_{j=1}^n{|Q(j)|^q} \Bigg)^{1/q} < \frac{|a_0|}{L} \, L|Q(0)| = |a_0Q(0)|\,.$$
Then $\sum_{j=0}^n{a_jQ(j)} \neq 0$, and hence Lemma 3.5 implies that 
$(x-1)^{m+1}$ does not divide $P$. We conclude that $\kappa_p(n,L) \leq m = \mu_q(n,L)$.
\qed \enddemo

\demo{Proof of Lemma 3.9}
Note that $\mu_q(n,K) \leq n$ for all $n \in {\Bbb N}$ and $L > 0$, as it is shown by $H \in {\Cal P}_n$ 
defined by $H(x) := \prod_{j=1}^n{(x-j)}$. 

\noindent Case 1: $0 < K < n^{-1/q}$. The choice $F \equiv 1$ gives the lemma. 

\noindent Case 2: $n^{-1/q} \leq K < 2$. Let $F$ be the polynomial given in Lemma 3.7 with $A := n$ and $M := (4K)^{2q}$. Then 
$$\split \sum_{j=1}^n{|F(j)|^q} & \leq \sum_{j \leq M^{-1/2}}{M^q} + \sum_{j > M^{1/2}}{\frac{1}{j^{2q}}} < 
M^{q-1/2} + \frac{1}{2q-1} \lfloor M^{-1/2} \rfloor^{-2q+1}  \cr 
& \leq (1 + 2^{2q-1})M^{q-1/2}\,, \cr \endsplit$$
so 
$$\Bigg( \sum_{j=1}^n{|F(j)|^q} \Bigg)^{1/q} < 4M^{1-1/(2q)} = K^{-1} F(0)\,,$$ 
and the degree $m$ of $F$ satisfies
$$m < \pi \sqrt{n} \root 4 \of M + 2 < \pi \sqrt{n}(4K)^{q/2} + 2 \leq \sqrt{n}(40K)^{q/2} + 2\,.$$

\noindent Case 3: $2 \leq K \leq \exp(n - 2\sqrt n)$. Let 
$R := \lfloor \sqrt{\log K} \rfloor + 1$, and let $F$ be the polynomial given in Lemma 3.7 with this $R$. Then 
$$|F(0)| > K \sum_{j=1}^n{|F(j)|} \geq K \Bigg( \sum_{j=1}^n{|F(j)|^q} \Bigg)^{1/q}\,,$$
and the degree $m$ of $F$ satisfies
$$m \leq \textstyle{\frac{44}{7}}R\sqrt{n} + 4 \leq 13 \sqrt{n\log K} + 4\,.$$ 

\noindent Case 4: $K > \exp(n - 2\sqrt n)$, $n \geq 9$. Then $\log K > n - 2\sqrt n \geq n/3$ for all 
$n \geq 9$. Hence the polynomial $F \in {\Cal P}_n$ defined by $F(x) := \prod_{j=1}^n{(x-j)}$ shows that
$$\mu_q(n,K) \leq n \leq \sqrt{3} \min \Big\{ \sqrt{n\log K},n \Big\}\,.$$  

\noindent Case 5: $K \geq 2$ and $n < 9$. Now the polynomial $F \in {\Cal P}_n$ defined by $F(x) := \prod_{j=1}^n{(x-j)}$ shows 
$$\mu_q(n,K) \leq n \leq 4\min \Big\{\sqrt{n\log K},n \Big\}\,.$$  
\qed \enddemo

\demo{Proof of Lemma 3.10}
First let $1 < K < 2$. Let $m = \lfloor \sqrt {n(K-1)/2} \rfloor + 1$.  Let $T_m$ be the Chebyshev polynomial of degree $m$ defined by
$$T_m(\cos t) = \cos(mt)\,, \quad t \in {\Bbb R}\,.$$
It is well known that $|T_m^\prime(1)| = m^2$ and $T_m^\prime(x)$ is increasing on $[1,\infty)$, hence 
$T_m(1+x) \geq 1 + m^2x$ for all $x > 0$. Now we define $F \in {\Cal P}_m$ by 
$$F(x) := T_m\left(\frac{-2x}{n-1}+ \frac{n+1}{n-1}\right) \,.$$   
Then $|F(x)| \leq 1$ for all $x \in [1,n]$, and 
$$F(0) \geq T_m\left(1+\frac{2}{n-1}\right) > 1 + \frac{m^2}{n-1} > 1 + \frac{m^2}{n} \geq K\,,$$ 
which finishes the proof in the case of $1 < K < 2$. Now let $k \geq 2$. Then the polynomial 
$F \in {\Cal P}_m$ chosen for $q = 1$, $n \in {\Bbb N}$, and $K \geq 2$ by Lemma 3.9 gives that 
$$|F(0)| > K \Bigg( \sum_{j=1}^n{|F(j)|^q} \Bigg)^{1/q} \geq K \max_{j \in \{1,2,\ldots, n\}}{|F(j)|}\,,$$
with  
$$m \leq 13 \min \Big \{\sqrt{n \log K}, n \Big\} + 4\,.$$
\qed \enddemo

\head 5. Proofs of the Theorems \endhead

\demo{Proof of Theorem 2.1}
Without loss of generality we may assume that $p \in (1,\infty)$, as the case $p = \infty$ follows by a 
simple limiting argument (or we may as well refer to the main result in [B-13]). By Lemma 3.6 we have  
$$\kappa_p(n,L) \leq \mu_q(n,L)$$
for every $n \in  {\Bbb N}$ and $L > 0$. The lower bounds for $\kappa_p(n,L)$ follows from Lemmas 3.2 and 3.3. 
The upper bounds for $\mu_q(n,L)$ follow from Lemma 3.9 with $K = L^{-1}$.
\qed \enddemo

\demo{Proof of Theorem 2.2}
By Lemma 3.6 we have
$$\kappa_1(n,L) \leq \mu_\infty(n,L)$$
for every $n \in  {\Bbb N}$ and $L > 0$.
The lower bounds for $\kappa_1(n,L)$ follow from Lemmas 3.3 and 3.4.
The upper bounds for $\mu_\infty(n,L)$ follow from Lemma 3.10 with $K = L^{-1}$.
\qed \enddemo

\enddocument